\documentclass[leqno]{macrorend}

\volumeyear{xxxx}\yearnumber{x}\volumenumber{xx}

\newtheorem{theo}{Theorem}[section]

\newtheorem{cor}{Corollary}[section]
\newtheorem{defi}{Definition}[section]

\newtheorem{exa}{Example}[section]

\newenvironment{sis}{\left\{\begin{aligned}}{\end{aligned}\right.}

\numberwithin{equation}{section}

\newcommand{\Z}{\mathbb{Z}}
\newcommand{\N}{\mathbb{N}}

\newcommand{\g}{\mathfrak{g}}

\newcommand{\U}{\mathfrak{U}}

\newcommand{\F}{\mathbb{F}}

\renewcommand{\O}{\mathcal{O}}

\usepackage{amsmath,amsthm,amsfonts,amssymb,amscd}
\input{xy}
\xyoption{all}

\begin{document} 

\selectlanguage{english}

\articolo[Simple finite group schemes and their infinitesimal deformations]{Simple finite group schemes and their infinitesimal deformations}{F. Viviani \footnotemark[1]}

\footnotetext[1]{Partially supported by FCT-Ci\^encia2008 (Portugal).}

\begin{abstract}
We show that the classification of simple finite group schemes over an algebraically closed 
field reduces to the classification
of abstract simple finite groups and of simple restricted Lie algebras in positive characteristic.
Both these two simple objects have been classified. We review this classification.
Finally, we address the problem of determining the infinitesimal deformations of simple
finite group schemes.  
\end{abstract}

%{\bf Mathematics Subject Classification:}  14L15, 17B70, 20D05, 17B50, 17B56. 

%{\bf Keywords:} Simple finite group schemes, simple abstract groups, simple Lie algebras,
%infinitesimal deformations.

%\maketitle

\section{Introduction}

In the first part of this paper, we show that a \emph{simple finite group scheme} over an 
algebraically closed field can be of two types:
either the constant group scheme associated to a simple (abstract) finite group or the
group scheme of height one associated to a simple 
restricted Lie algebra. The classification of these two kind of simple objects 
(the simple Lie algebras have been classified only for $p\neq 2,3$)
was certainly among the greatest achievements of the mathematics of the last 
century.

% there was a complete classification of these two kind of simple objects
%(the simple Lie algebras have been classified only for $p\neq 2,3$).

\emph{Simple finite groups} were classified during the years 1955-1985
thanks to the contribution of many mathematicians 
(see \cite{SOL1}, \cite{SOL2}, \cite{ASC} for a nice historical account).

\emph{Simple Lie algebras} over an algebraically closed field 
$F$ of characteristic $p\neq 2,3$ have recently been classified by Block-Wilson-Premet-Strade 
(see  \cite{BW1}, \cite{SW}, \cite{STR}).
The classification says that for $p\geq 7$ the simple Lie algebras
can be of two types: of classical type and of generalized Cartan type. 
The algebras of \emph{classical type} are obtained by considering the simple Lie 
algebras in characteristic zero, by taking a model 
over the integers and then reducing modulo the prime $p$.
The algebras of \emph{generalized Cartan type} are the finite-dimensional 
analogues of the four classes of infinite-dimensional complex simple Lie algebras, 
which occurred in Cartan's classification of Lie pseudogroups.
In characteristic $p=5$, apart from the above two types of algebras, there is one 
more family of simple Lie algebras called Melikian algebras.
In characteristic $p=2, 3$, there are many exceptional simple restricted Lie algebras 
and the classification seems still far away.

After passing in review these two classification results, in the last part of this
paper we address the problem of determining the \emph{infinitesimal 
deformations} of such simple finite group schemes. The group schemes 
associated to simple finite (abstract) groups and to simple Lie 
algebras of classical type are known to be rigid apart from some bad characteristic 
of the base field.
We show that this is not the case for the group schemes associated to simple 
Lie algebras of Cartan type. In particular we 
determine the infinitesimal deformations of the group schemes of height one
associated to the restricted simple Lie algebras of Cartan type. 
We do this by computing the second restricted cohomology group
of these algebras with values in the adjoint representation.
It remains an open problem to extend the above results to the other simple
restricted Lie algebras.

\section{Finite group schemes}

Let $k$ be a field of characteristic $p\geq 0$.
There are three equivalent ways to define a finite group scheme over $k$.

\begin{defi}
A finite group scheme $G$ over $k$ can be defined equivalently as
\begin{itemize}
\item[(i)] A finite scheme $G$ which is a group object in the category of schemes;
\item[(ii)] A scheme of the form $G={\rm Spec}(A)$, where $A$ is a finite dimensional
commutative $k$-Hopf algebra;  
\item[(iii)] A scheme $G$ whose functor of points $F_G$ takes values in finite groups, that is
$$\begin{aligned}
F_G:\{k-{\rm schemes}\}&\to \{{\rm finite \:\:groups}\}\\
S& \mapsto G(S):={\rm Hom}_k(S, G).
\end{aligned} 
$$
\end{itemize}
\end{defi}
The following Theorem collects the basic structure of finite group schemes.

\begin{theo}\label{structure}
Let $G$ be a finite group scheme over the field $k$. 
\begin{itemize}
\item[(1)] If ${\rm char}(k)=0$ then $G$ is \'etale.
\item[(2)] If ${\rm char}(k)=p>0$ then there is a unique extension 
$$1\to G^0\to G \to G^{et}\to 1$$ 
such that $G^0$ is connected and $G^{et}$ is \'etale. If $k$ is perfect then the above exact 
sequence splits, that is $G=G^0\rtimes G^{et}$.
\item[(3)] If ${\rm char}(k)=p>0$ and $G$ is connected, there is a minimal natural number $n\geq 1$
such that $G={\rm Ker}(F^n)$ (it is called the height of $G$), where $F^n:G \to G^{(p^n)}$
is the the $n$-th iteration of the Frobenius map.
In the filtration of normal closed subgroups:
$$1\lhd {\rm Ker}(F) \lhd {\rm Ker}(F^2)\lhd \cdots \lhd {\rm Ker}(F^n)=G^0,$$ 
each factor ${\rm Ker}(F^{i+1})/{\rm Ker}(F^i)$ has height one (i.e. it 
has vanishing Frobenius).
\end{itemize}
\end{theo} 
\begin{proof}
\begin{itemize}
\item[(1)] By a result of Cartier, every $k$-group scheme is smooth if ${\rm char}(k)$ $=0$
(see \cite[Chap. 11.4]{WAT}). Therefore every finite group scheme in characteristic 
zero is \'etale.
\item[(2)] The above exact sequence is obtained by taking $G^0$ the connected component of $G$
containing the identity and $G^{et}={\rm Spec}(\pi_0(k[G]))$, where $\pi_0(k[G])$ is the maximal
separable sub-Hopf algebra of the algebra $k[G]$ of regular functions on $G$.
If $k$ is perfect, then $G^{et}\cong {\rm Spec}(k[G]_{red})$, where
$k[G]_{red}$ is the maximal reduced quotient of $k[G]$, which gives the splitting 
(see  \cite[Chap. 6]{WAT}).
\item[(3)] The kernel of $F^n$ is represented by the Hopf algebra $k[G]/(x^{p^n}\: | \: x\in I)$,
where $I$ is the augmentation ideal of $k[G]$ (that is, the maximal ideal corresponding 
to the origin of $G$). The first assertion follows from the fact that if $G$ is connected then
the augmentation ideal $I$ is nilpotent. The second assertion is clear.
\end{itemize}
\end{proof}

Therefore every finite group scheme can be realized as extensions of \'etale and height one 
finite group schemes (these latter occur only if ${\rm char}(k)=p>0$). 
Now we provide an explicit description of these two building blocks.

\begin{theo}\label{etale-gs}
If $k$ is algebraically closed the we have a bijection
$$\{\text{\'Etale $k$-group schemes}\}\longleftrightarrow \{\text{Finite (abstract) groups}\}$$ 
\end{theo}
\begin{proof}
The bijection is realized explicitly as follows: to an \'etale $k$-group scheme $G$ one associates 
the finite group $G(k)$ of its $k$-points. 
Conversely, to an abstract finite group $\Gamma$ one 
associates the finite group scheme whose $k$-Hopf algebra is the $k$-algebra $k^{\Gamma}$ of  functions from 
$\Gamma$ to $k$ with comultiplication given by $\Delta(e_{\rho})=\oplus_{\rho=\sigma\cdot \tau}
e_{\sigma}\otimes e_{\tau}$, where $e_{\rho}$ is the function sending $\rho\in \Gamma$ into 
$1$ and the other elements of $\Gamma$ to $0$. One can check the above maps are one the 
inverse of the other (see \cite[Chapter 6.4]{WAT}).
\end{proof}

In order to describe the finite group schemes of height one, we need to recall the definition
of the restricted Lie algebras (sometimes called $p$-Lie algebras) over a field $k$ of positive
characteristic.

\begin{defi}[\cite{JAC2}]
A Lie algebra $L$ over a field $k$ of characteristic
$p>0$ is said to be \emph{restricted} (or a $p$-Lie algebra) if it is endowed
with a map (called $p$-map) $[p]:L\to L$, $x\mapsto x^{[p]}$, 
which satisfies the following conditions:
\begin{enumerate}
\item ${\rm ad}(x^{[p]})={\rm ad}(x)^{[p]}$ for every $x\in L$.
\item $(\alpha x)^{[p]}=\alpha^p x^{[p]}$ for every $x\in L$ and every $\alpha\in k$.
\item $(x_0+x_1)^{[p]}=x_0^{[p]}+x_1^{[p]}+\sum_{i=1}^{p-1} 
s_i(x_0,x_1)$ for every $x,y\in L$,
where the element $s_i(x_0,x_1)\in L$ is defined by
$$s_i(x_0,x_1)=-\frac{1}{r}\sum_u {\rm ad}x_{u(1)}\circ {\rm ad}x_{u(2)}\circ
\cdots \circ {\rm ad}x_{u(p-1)}(x_1),$$
the summation being over all the maps $u:[1,\cdots,p-1]\to \{0,1\}$ taking $r$-times
the value $0$.
\end{enumerate}
\end{defi}

The last two conditions in the above definition can be rephrased by saying that the map
$ a \mapsto a^p-a^{[p]}$ from $L$ into the universal enveloping algebra $\U_L$ is p-semilinear, 
where $a^p$ denotes the $p$-th self product of $a$ in $\U_L$. 
We give examples of restricted Lie algebras. 

\begin{exa}
\noindent 
\begin{enumerate}

\item Let $A$ an associative $k$-algebra, where ${\rm char}(k)=p>0$. 
Then the Lie algebra ${\rm Der}_F A$ of 
$k$-derivations of $A$ into itself is a restricted Lie algebra 
with bracket $[D_1,D_2]=D_1\circ D_2-D_2\circ D_1$ and $p$-map
$D\mapsto D^p:=D\circ \cdots\circ D$.

\item Let $G$ a group scheme over $k$, where ${\rm char}(k)=p>0$.
Then the Lie algebra ${\rm Lie}(G)$ associated 
to $G$ is a restricted Lie algebra with respect to the $p$-map given by the differential 
of the homomorphism $G\to G$, $x\mapsto x^p:=x\circ \cdots\circ x$.

\end{enumerate}
\end{exa}

The Lie algebra ${\rm Lie}(G)$ associated to a $k$-group scheme $G$ in positive 
characteristic $p$ (as in the above example) depends only on the first height truncation 
${\rm Ker}(F)\lhd G$ of $G$.  Indeed, we have the following

\begin{theo}\label{heightone-gs}
Let $k$ be a field of characteristic $p>0$. Then there is a bijection 
$$
\mbox{\{Restricted $k$-Lie algebras\}}  \longleftrightarrow
\mbox{\{Finite $k$-group schemes of height one\}}
$$
\end{theo}
\begin{proof}
The bijection is realized explicitly as follows: 
to a finite group scheme $G$ of height $1$, one associates 
the restricted Lie algebra ${\rm Lie}(G)$ with $p$-map given by the differential
of the map $x\mapsto x^p$. Conversely, to a restricted Lie
algebra $(L,[p])$, one associates the finite group scheme corresponding to the dual
of the restricted enveloping Hopf algebra $\U^{[p]}(L):=\U(L)/(x^p-x^{[p]})$. 
One can check that the above map are inverse one of the other 
(see \cite[Chapter 2.7]{DG}).
\end{proof}

As a consequence of the structure theorems for finite group schemes discussed in this section, we 
get

\begin{cor}
A simple finite group scheme over an algebraically closed field $k$ is either the 
\'etale group scheme associated to a simple (abstract) finite group or 
(if ${\rm char}(k)=p>0$) the height one group scheme associated to a simple restricted 
Lie algebra. 
\end{cor}

Both these two simple objects have been classified (the last ones assuming that 
$p\neq 2, 3$)!

\section{Classification of simple finite groups}

The simple finite groups are divided in four families:

\begin{itemize}
 \item  Cyclic groups $\Z/p\Z$ with prime order; 
\item Alternating groups $A_n$, $n\geq 5$;
\item  Simple groups of Lie type, including 

 \begin{itemize}
  \item Classical groups of Lie type:
 
  \begin{itemize}
\item[(1)] $A_n(q)$: The Projective Special Linear Group ${\rm PSL}_{n+1}(\F_q)$.
\item[(2)] ${}^2A_n(q^2)$: The Projective Special Unitary Group ${\rm PSU}_{n+1}(\F_{q^2})$,
 with respect to the Hermitian form on $\F_{q^2}^{n+1}$
$$\Psi(w, v)=\sum_{i=1}^{n+1} w_i^q v_i,$$ 
where $w=(w_1, \cdots, w_{n+1})$ and $v=(v_1, \cdots, v_{n+1})$.
\item[(3)] $B_n(q)$:  The subgroup $O_{2n+1}(\F_q)$ of the special orthogonal group in dimension 
$2n+1$ formed by the elements having spinor norm $1$.
\item[(4)] $C_n(q)$: The projective Symplectic group ${\rm PSp}_{2n}(\F_q)$.
\item[(5)] $D_n(q)$: The subgroup $O^+_{2n}(\F_q)$ of the projective special split 
orthogonal group in dimension $2n$ formed by the elements having spinor norm $1$.
\item[(6)] ${}^2D_n(q)$: The subgroup $O^-_{2n}(\F_{q^2})$ of elements of spinor norm $1$ 
in the projective special non-split orthogonal group in dimension $2n$. 
\end{itemize}
    
%  The classical groups, namely the groups of projective special linear, unitary, symplectic, or  orthogonal transformations over a finite field;
 
 \item Exceptional and twisted groups of Lie type, obtained via the Steinberg construction,
starting with an automorphism of a Dynkin diagram and an automorphism of a finite field.
The resulting groups are called: ${}^2B_2(2^{2n+1})$, ${}^3D_4(q^3)$, $E_6(q)$, 
${}^2E_6(q^2)$, $E_7(q)$, $E_8(q)$, $F_4(q)$, ${}^2F_4(2^{2n+1})$, $G_2(q)$, 
${}^2G_2(3^{2n+1})$ (see \cite{CA}).

\end{itemize}
\item $26$ sporadic groups (see \cite{WIL}).
\end{itemize}

%The {\bf classical groups} of Lie type are

\section{Classification of simple restricted Lie algebras}

First of all, we show that there is a bijection between simple restricted Lie algebras,
that is restricted Lie algebras without restricted ideals 
(i. e.  ideals closed under the $p$-map),
and simple Lie algebras (not necessarily restricted).

\begin{theo}
There is a bijection
$$\{\text{Simple restricted Lie algebras}\}\longleftrightarrow \{\text{Simple Lie algebras}\}.$$
Explicitly to a simple restricted Lie algebra $(L,[p])$ we associate its derived algebra
$[L,L]$. Conversely to a simple Lie algebra $M$ we associate the restricted subalgebra 
$M^{[p]}$ of 
${\rm Der}_F(M)$ generated by ${\rm ad}(M)$ (which is called the universal $p$-envelope of 
$M$).
\end{theo}
\begin{proof}
We have to prove that the above maps are well-defined and are inverse one of the other.

$\bullet$ Consider a simple restricted Lie algebra $(L,[p])$. The derived subalgebra 
$[L,L]\lhd L$ is a non-zero ideal (since $L$ is not abelian) and therefore 
$[L,L]_p=L$, where $[L,L]_p$ denotes the $p$-closure of $[L,L]$ inside $L$.

Take a non-zero ideal $0\neq I\lhd [L,L]$. 
Since $[L,L]_p=L$, we deduce from \cite[Chapter 2, Prop. 1.3]{FS} that $I$ is also an 
ideal of $L$ and therefore $I_p=L$ by the restricted simplicity of $(L,[p])$.
From loc. cit., it follows also that $[L,L]=[I_p,I_p]=[I,I]\subset I$ from which 
we deduce that $I=L$. Therefore $[L,L]$ is simple. 

Since ${\rm ad}:L\to {\rm Der}_F(L)$ is injective and $[L,L]_p=L$, it follows 
by loc. cit. that ${\rm ad}:L\to {\rm Der}_F([L,L])$ is injective. Therefore we have that
$[L,L] \subset L \subset {\rm Der}_F([L,L])$ and hence $[L,L]^{[p]}=[L,L]_p=L$.

$\bullet$ Conversely, start with a simple Lie algebra $M$ and consider its universal
$p$-envelop $M<M^{[p]}<{\rm Der}_F(M)$. 

Take any restricted ideal $I\lhd_p M^{[p]}$. By loc. cit., 
we deduce $[I,M^{[p]}]\subset I\cap [M^{[p]}, M^{[p]}]=I\cap [M,M]=I\cap M\lhd M$. 
Therefore, by the simplicity of $M$, either $I\cap M=M$ or $I\cap M=0$. 
In the first case, we have that $M\subset I$ and therefore $M^{[p]}=I$.
In the second case, we have that $[I,M^{[p]}]=0$ and therefore $I=0$
because $M^{[p]}$ has trivial center. We conclude that $M^{[p]}$ is simple restricted.

Moreover, by loc. cit., we have that $[M^{[p]},M^{[p]}]=[M,M]=M$.
\end{proof}
Note that the intersection of the two types of Lie algebras appearing in the above 
correspondence is the set of restricted simple Lie algebras, namely those restricted 
Lie algebra which does not have any proper ideal.

Simple Lie algebras (and their minimal $p$-envelopes) over $k=\overline{k}$ of ${\rm char}(k)=p>3$ 
have been classified by Block-Wilson-Premet-Strade (1984-2005), answering to
a conjecture of Kostrikin-Shafarevich (1966).
They are divided into two types:
\begin{itemize}
\item Lie algebras of \emph{Classical type}; 
\item Lie algebras of \emph{(generalized) Cartan type}.
\end{itemize}

\subsection{Lie algebras of classical type}

The Lie algebras of classical type are reduction of simple
Lie algebras in characteristic zero.
Simple Lie algebras over an algebraically closed field of {\bf characteristic zero} were
classified at the beginning of the XIX century by Killing and Cartan.
The classification proceeds as follows:
first the non-degeneracy of the Killing form is used to establish a correspondence between 
simple Lie algebras and irreducible root systems and  then the irreducible root systems  are 
classified by mean of their associated Dynkin diagrams. 
Explicitly:
$$\begin{aligned} \text{DYNKIN DIAGRAMS} & \longleftrightarrow  \text{SIMPLE LIE ALGEBRAS}\\
A_n \: (n\geq 1) & \hspace{1cm} \mathfrak{sl}(n+1)  \\
B_n \: (n\geq 2) & \hspace{1cm} \mathfrak{so}(2n+1)  \\
C_n \: (n\geq 3) & \hspace{1cm} \mathfrak{sp}(2n)  \\
D_n \: (n\geq 4) & \hspace{1cm} \mathfrak{so}(2n)  \\
E_6, E_7, E_8, F_4, G_2 & \hspace{1cm} \mbox{Exceptional Lie algebras}\\
\end{aligned}$$
where $\mathfrak{sl}(n+1)$ is the special linear algebra, $\mathfrak{so}(2n+1)$ is 
the special orthogonal algebra of odd rank, $\mathfrak{sp}(2n)$ is the symplectic algebra and  
$\mathfrak{so}(2n)$ is the special orthogonal algebra of even rank.
For the simple Lie algebras corresponding to the 
exceptional Dynkin diagrams, see the book \cite{JACEXE} or the nice account in \cite{BAE}.

These simple Lie algebras admit a model over the integers via the (so-called) Chevalley bases.
Therefore, via reduction modulo a prime $p$, one obtains a restricted 
Lie algebra over $\mathbb{F}_p$, which is simple up to a quotient by a small ideal. 
For example $\mathfrak{sl}(n)$ is not simple 
if $p$ divide $n$, but its quotient $\mathfrak{psl}(n)=\mathfrak{sl}(n)/(I_n)$ by the unit 
matrix $I_n$ becomes simple.
There are similar phenomena occurring only for $p=2, 3$ for the other Lie 
algebras (see \cite[Page 209]{STR}). The restricted simple algebras obtained in this way are 
called algebras of {\bf classical type}.  
Their Killing form is non-degenerate except at a finite number of primes.
Moreover, they can be characterized 
as those restricted simple Lie algebras admitting a projective representation with nondegenerate 
trace form (see \cite{BLO}).

\subsection{Lie algebras of Cartan type}

However, there are restricted simple Lie algebras which have
no analogous in characteristic zero and therefore are called non-classical.
The first example of a non-classical restricted simple Lie algebra is due to E. Witt, who in 1937 
realized that the derivation algebra $W(1):= {\rm Der}_k(k[X]/(X^p))$ over a field $k$ of 
characteristic $p>3$ is simple with a degenerate Killing form.
In the succeeding three decades, many more non-classical restricted simple Lie algebras 
have been found (see \cite{JAC1}, \cite{FRA1}, \cite{AF}, \cite{FRA2}). 
The first comprehensive conceptual approach to constructing these 
non-classical restricted simple Lie algebras was proposed by Kostrikin-Shafarevich 
and Kac (see \cite{KS}, \cite{KS2}, \cite{KAC1}). 
They showed that all the known examples  can be 
constructed as finite-dimensional analogues of the four 
classes of infinite-dimensional complex simple Lie algebras, which occurred in Cartan's 
classification of Lie pseudogroups (see \cite{CAR}).

Denote with $\O(m)$ the divided power $k$-algebra in $m$-variables. It is the commutative and 
associative algebra with unit defined by the generators $x^{\alpha}$ for $\alpha\in \N^m$
satisfying the relations 
$$ x^{\alpha}\cdot x^{\beta}=\binom{\alpha+\beta}{\alpha} x^{\alpha+\beta}
:=\prod_{i=1}^m \binom{\alpha_i+\beta_i}{\alpha_i} x^{\alpha+\beta}.$$

\noindent For any $m$-tuple $\underline{n}\in \N^m$, we define the truncated subalgebra
of $\O(m)$ 
$$\O(m;\underline{n})={\rm span}
\langle x^{\alpha}\, | \, 0 \leq \alpha_i < p^{n_i}\rangle.$$

The simple Lie algebras of Cartan type (over an algebraically closed field $k$ of
characteristic $p>3$) are divided in four families, called Witt-Jacobson, Special,
Hamiltonian and Contact algebras, plus an exceptional family of Lie algebras
in characteristic $p=5$, called Melikian algebras. We list the simple \emph{graded}
Lie algebras of Cartan type. The general simple Lie algebras of Cartan type are filtered
deformations of these graded Lie algebras (see \cite[Chap. 6]{STR}).

\begin{itemize}

\item[(1)] WITT-JACOBSON: The Witt-Jacobson Lie algebra $W(m;\underline{n})$
is the subalgebra of ${\rm Der}_k \O(m;\underline{n})$ of special
derivations:

$$W(m;\underline{n})=\{D\in {\rm Der}_k \O(m;\underline{n})\: : D(x^{(a)})= \sum_{i=1}^m 
x^{(a-\epsilon_i)}D(x_i)\}. $$
The Witt-Jacobson Lie algebra is a free $\O(m, \underline{n})$ generated by the 
special derivations $\partial_i$ defined by $\partial_i(x^{(a)})=x^{(a-\epsilon_i)}$.

\item[(2)] SPECIAL: The Special Lie algebra is the subalgebra of $W(m, \underline{n})$,
$m\geq 3$, 
of derivations preserving the special form $\omega_S= dx_1\wedge\cdots\wedge dx_m$:
$$S(m;\underline{n})^{(1)}=\{D\in W(m;\underline{n})\: | \: D(\omega_S)=0\}^{(1)},$$
where $(1)$ denotes the derived algebra.

 \item[(3)] HAMILTONIAN: The Hamiltonian algebra is the subalgebra of $W(2r,\underline{n})$ 
of derivations preserving the Hamiltonian form 
$\omega_H=dx_1\wedge dx_{r+1}+\cdots +dx_r\wedge dx_{2r}$:
$$H(2r;\underline{n})^{(2)}=\{D\in W(2r;\underline{n})\: | \: D(\omega_H)=0\}^{(2)},$$
where $(2)$ denotes the double derived algebra.

 \item[(4)] CONTACT: The Contact algebra is the subalgebra of $W(2r+1,\underline{n})$ 
of derivations preserving the contact form 
$\omega_K=\sum_{i=1}^r(x_i dx_{i+r}-x_{i+r}dx_i)+dx_{2r+1}$ up to multiples:
$$K(2r+1;\underline{n})^{(1)}=\{D\in W(2r+1;\underline{n})\: | \:
D(\omega_K)\in \O(2r+1;\underline{n}) \omega_K \}^{(1)}.$$

\item[(5)] MELIKIAN (only if ${\rm char}(k)=5$): The Melikian algebra (introduced in \cite{MEL}) is defined as 
$$M(n_1,n_2)=\O(2, (n_1, n_2))\oplus W(2, (n_1, n_2))\oplus \widetilde{W(2, (n_1, n_2))},$$
where $\widetilde{W(2, (n_1, n_2))}$ denote a copy of $W(2, (n_1, n_2))$ and the Lie bracket 
is defined by the following rules (for all $D, E\in W(2, (n_1, n_2))$ and 
$f, g\in \O(2, (n_1, n_2))$):
$$\begin{sis}
& [D, E]:= [D; E], \\
&[D,\widetilde{E}]:=\widetilde{[D,E]}+2\,{\rm div}(D)\widetilde{E},\\
&[D,f]:=D(f)-2\, {\rm div}(D)f,\\
&[f_1\widetilde{D_1}+f_2\widetilde{D_2},g_1\widetilde{D_1}+g_2\widetilde{D_2}]:=
f_1g_2-f_2g_1,\\
&[f,\widetilde{E}]:=f E,\\
&[f,g]:=2\,(gD_2(f)-fD_2(g))\widetilde{D_1}+2\,(fD_1(g)-gD_1(f))\widetilde{D_2},\\
\end{sis}$$
where ${\rm div}(f_1D_1+f_2D_2):=D_1(f_1)+D_2(f_2)\in \O(2, (n_1, n_2))$.

\end{itemize}

The above algebras are restricted if and only if $\underline{n}=\underline{1}$.
In general, denoting with $X(m, \underline{n})$ one of the above simple Lie algebras,
its minimal $p$-envelope $X(m, \underline{n})_{[p]}$ inside its derivation 
algebra ${\rm Der}_k(X(m, \underline{n}))$ is equal to (see \cite[Theo. 7.2.2]{STR}):

$$X(m;\underline{n})_{[p]}=X(m;\underline{n})+\sum_{i=1}^m \sum_{0< j_i <n_i} k \cdot 
\partial_i^{p^{j_i}}.$$ 

In this way we obtain all the simple restricted Lie algebras,
up to filtered deformations. We observe that in each of the above simple Lie 
algebras of Cartan type, the Killing form is always degenerate.

\section{Deformations of simple finite group schemes}

Having at our hand a complete classification of the simple finite group schemes, it is natural  to study 
their properties. Here we are interested in their deformations.
%more deeply the properties of the simple finite group schemes. 
%One of these properties is to ask about their deformations. 
According to Grothendieck's philosophy,
one should first understand the infinitesimal deformations.
It is a classical result (of difficult attribution) that the infinitesimal
deformations of the finite group scheme $G$ are parametrized by $H^2(G, \g)$,
the second cohomology group of $G$ with values in the adjoint representation
$\g={\rm Lie}(G)$.

\subsection{Deformations of simple (abstract) groups}

From the Maschke's theorem, it follows that every constant group scheme is rigid
over a field of characteristic not dividing the order of the group. In particular this applies
to the finite group schemes associated to simple finite groups.
We don't known what happens if the characteristic of $k$ divides the order of the group.

\subsection{Deformations of Lie algebras of Classical type}

If $G$ is a finite group scheme of height one associated to a restricted Lie algebra
$(\g, [p])$, the cohomology $H^i(G; \g)$ of $G$ with values in the adjoint representation 
is equal to the restricted cohomology $H_*^i(\g, \g)$ of $\g$ with values in the adjoint 
representation (see \cite{HOC}).

The restricted cohomology is related to the ordinary cohomology by two spectral sequences 
(see \cite{JAN} and \cite{FAR}):
$$\begin{sis}
&E_1^{p,q}={\rm Hom}_{Fr}(S^p \g, H^{q-p}(\g,M))\Rightarrow H_{*}^{p+q}(\g,M) \:
\text{ if } p\neq 2,  \hspace{0,5cm} & \\
&E_2^{p,q}={\rm Hom}_{Fr}(\Lambda^q \g, H_{*}^{p}(\g,M))\Rightarrow H^{p+q}(\g,M),
\hspace{1cm} & 
\end{sis}$$
where $S^p\g$ and $\Lambda^q\g$ denote, respectively, the $p$-th symmetric power and 
the $q$-th alternating power, and ${\rm Hom}_{Fr}$ denotes the homomorphisms which are 
Frobenius linear. In the case of the adjoint representation and for centerless Lie algebras,
the above spectral sequences give rise to the following relations for the low
cohomology groups:
$$\begin{sis}
&  H_*^1(\g,\g)=H^1(\g,\g), \\
&  0 \to H_*^2(\g, \g)\to H^2(\g,\g)\to {\rm Hom}_{\rm Frob}(\g,H^1(\g,\g)).
\end{sis}$$

\noindent It is a classical fact that simple Lie algebras are rigid in characteristic zero.

\begin{theo}[Whithead] If $\g$ is a simple Lie algebra over a field $k$ of 
${\rm char}(k)=0$, then $H^*(\g,\g)=0$. In particular $\g$ is rigid. 
\end{theo}

The proof uses the non-degeneracy of the Killing form.
If ${\rm char}(k)=p$ does not divide the discriminant of the Killing form, then the same proof
gives the rigidity of simple Lie algebras of classical type in characteristic $p$. 
Indeed Rudakov (\cite{RUD}) has shown the following 
 
\begin{theo}[Rudakov] If $k$ is a field of ${\rm char}(k)=p\geq 5$ and $\g$ is a simple
Lie algebra of classical type, then $\g$ is rigid.
\end{theo}
 
Note, however, that the preceding result is false if $p=2,3$ (see  \cite{Che2},
\cite{Che3}, \cite{Che1}).

\subsection{Deformations of Lie algebras of Cartan type}

The Lie algebras of Cartan type do have infinitesimal deformations, unlike the Lie algebras of 
classical type. We computed in \cite{VIV5} (building upon \cite{VIV1},
\cite{VIV2} and \cite{VIV3}) the infinitesimal deformations
of the simple finite group schemes associated to simple and restricted Lie algebras
of Cartan type.

\begin{theo}[\cite{VIV5}]
$$\begin{sis}
& h_*^2(W(m;\underline{1}),W(m;\underline{1}))=m, \\
& h_*^2(S(m;\underline{1}),S(m;\underline{1}))=m, \\
& h_*^2(H(2r;\underline{1}),H(2r;\underline{1}))=2r+1, \\
& h_*^2(K(2r+1;\underline{1}),K(2r+1;\underline{1}))=2r+1, \\
& h_*^2(M(1,1),M(1, 1))=5.
\end{sis}$$
\end{theo}

Moreover, in each of the above cases, we get explicit generators for the above 
second cohomology groups.
It would be interesting to extend the above computation to the others simple
restricted Lie algebras of Cartan type (i.e. the ones with $\underline{n}\neq \underline{1}$).

\subsection{Further developments}

We hope that, using similar techniques to the ones used in \cite{VIV5},
it could be possible to answer to the following Questions (for all simple finite group schemes $G$):

\begin{itemize}

 \item What is the space of \emph{obstructions} $H^3(G,\g)$?

 \item What is the  \emph{semi-versal deformation space} of $G$?
 %that is going from infinitesimal to global formal deformations.

\end{itemize}

%\begin{thebibliography}{99} 

%\bibitem{1}{D.H. Ackeley, G.E. Hilton and T.J. Sejnovski, A learning algorithm for Bolzmann machine, \em 
%Cognitive Science, }{\bf 62} (1985), 147 - 169. 

%\bibitem{2}{ D.O. Hebb, \em Organization of Behaviour,} Wiley, New York, 1949. 

%\end{thebibliography} 

\bigskip

\begin{flushleft}

{\bf AMS Subject Classification: 14L15, 17B70, 20D05, 17B50, 17B56}\\[2ex]

% Write more than one author separately if they have different 
% affiliations, otherwise write the names on the same line, separeted 
% by commas.
%
Filippo Viviani,\\
Dipartimento di Matematica\\
Universit\`a Roma Tre,
Largo S. Leonardo Murialdo 1,
00146 Roma (Italy) \\
e-mail: \textit{viviani@mat.uniroma3.it}\\[2ex]

% leave it blanck, you dont know these infos yet
%
\textit{Lavoro pervenuto in redazione il MM.GG.AAAA.}

\end{flushleft}

\normalsize
\label{\thechapter:lastpage}


\footnotesize


\begin{thebibliography}{INTRO}

\addcontentsline{toc}{section}{References}

%\bibitem{1}{D.H. Ackeley, G.E. Hilton and T.J. Sejnovski, A learning algorithm for Bolzmann machine, \em 
%Cognitive Science, }{\bf 62} (1985), 147 - 169. 

\bibitem{AF}{A. A. Albert and M. S. Frank, Simple Lie algebras of 
characteristic $p$, \em Rend. Sem. Mat. Univ. Politec. Torino, } {\bf 14} (1954-1955), 117--139.

\bibitem{ASC}{M. Aschbacher, The status of the classification of the finite 
simple groups \em Notices Amer. Math. Soc.,} {\bf 51}  (2004), 736--740.


\bibitem{BAE}{ J. Baez, The Octonions, \em Bull. Amer. Math. Soc., } {\bf 39} (2002), 
145-205.

\bibitem{BLO}{ R. E. Block, Trace forms on Lie algebras, \em Canad. J. Math.,}
{\bf 14} (1962), 553--564.

\bibitem{BW1}{ R. E. Block and R. L. Wilson, The restricted simple Lie algebras
are of classical or Cartan type, \em Proc. Nat. Aca. Sci. USA,} {\bf 81} (1984), 
5271--5274.

%\bibitem[BW88]{BW2} R.E. Block and R.L. Wilson: \emph{Classification of the restricted %simple Lie algebras}.  J. Algebra  114  (1988), 115--259.

\bibitem{CAR}{ E. Cartan, Les groupes de transformations continus, 
 infinis, simples, \em Ann. Sci. \'Ecole Nor. Sup.,} {\bf 26} (1909), 93--161.

\bibitem{CA}{R. W. Carter, \em Simple groups of Lie type,} 
Pure and Applied Mathematics, Vol. {\bf 28}, John Wiley and Sons, 
London-New York-Sydney, 1972.

\bibitem{Che1}{ N. G. Chebochko, Deformations of classical Lie algebras with a
 homogeneous root system in characteristic two I (Russian), \em
Mat. Sb., } {\bf 196,}  (2005), 125--156. English translation:  {\em Sb. Math.} {\bf 196 }
(2005), 1371--1402.

\bibitem{Che2}{ N. G. Chebochko and M. I. Kuznetsov, Deformations of classical Lie
algebras (Russian),\em Mat. Sb.,} {\bf 191} (2000), 69--88.
English translation: {\emph Sb. Math.} {\bf 191} (2000), 1171--1190.

\bibitem{Che3}{ N. G. Chebochko, S. A. Kirillov and M. I. Kuznetsov,
Deformations of a Lie algebra of type $G\sb 2$ of characteristic three (Russian),
\em Izv. Vyssh. Uchebn. Zaved. Mat.,} 2000, 33--38.
English translation: {\em Russian Math. (Iz. VUZ)} {\bf 44} (2000), 31--36.


%\bibitem[CEL70]{CEL} M. Ju. Celousov: \emph{Derivations of Lie algebras of
% Cartan type} (Russian).
% Izv. Vys\v s. U\v cebn. Zaved. Matematika 98 (1970), 126--134.


%\bibitem[CE48]{CE} C. Chevalley and S. Eilenberg: \emph{Cohomology theory of Lie groups
%and Lie algebras}. Trans. Amer. Math. Soc. 63 (1948), 85--124.

\bibitem{DG}{ M. Demazure and P. Gabriel: \em Groupes alg\'ebriques,
Tome I: G\'eom\'etrie alg\'ebrique, g\'en\'eralit\'es, groupes commutatifs} (French),
Avec un appendice  Corps de classes local par Michiel Hazewinkel,
Masson and Cie, Editeur, Paris; North-Holland Publishing Co., Amsterdam, 1970.

%\bibitem[DK78]{DK} A. S. D\v zumadildaev and A. I. Kostrikin:
%\emph{Deformations of the Lie algebra $W\sb{1}(m)$} (Russian).
%Algebra, number theory and their applications.  Trudy Mat. Inst. Steklov.  148  (1978),
%141--155, 275.

%\bibitem[DZU80]{DZ1} A. S. D\v zumadildaev: \emph{Deformations of general Lie algebras
%of Cartan type} (Russian).  Dokl. Akad. Nauk SSSR  251  (1980), no. 6, 1289--1292.
%English translation: Soviet Math. Dokl. 21 (1980), no. 2, 605--609.

%\bibitem[DZU81]{DZ2} A. S. D\v zumadildaev: \emph{Relative cohomology and deformations
%of the Lie algebras of Cartan types} (Russian).  Dokl. Akad. Nauk SSSR  257  (1981), no. %5, 1044--1048.
%English translation: Soviet Math. Dokl. 23 (1981), no. 2, 398--402.

%\bibitem[DZU89]{DZ3} A. S: D\v zumadildaev: \emph{Deformations of the Lie algebras
%$W\sb n(m)$} (Russian).  Mat. Sb.  180  (1989),  no. 2, 168--186.
%English translation: Math. USSR-Sb.  66  (1990),  no. 1, 169--187.


\bibitem{FAR}{ R. Farnsteiner, Cohomology groups of reduced enveloping algebras,
\em Math. Zeit.,} {\bf 206}  (1991), 103--117. 

\bibitem{FS}{ R. Farnsteiner and H. Strade, \em Modular Lie algebras
and their representation,} Monographs and textbooks in pure and applied mathematics,
vol. {\bf 116}, Dekker, New York, 1988.

\bibitem{FRA1}{ M. S. Frank, A new class of simple Lie algebras, \em
Proc. Nat. Acad. Sci. USA,} {\bf 40} (1954), 713--719. 

\bibitem{FRA2}{ M. S. Frank, Two new classes of simple Lie algebras,
\em Trans. Amer. Math. Soc.,} {\bf 112} (1964), 456--482. 

%\bibitem[GER64]{GER1} M. Gerstenhaber: \emph{On the deformation of rings and algebras}.
%Ann. of Math. 79  (1964), 59--103.

% \bibitem[GER66]{GER2} M. Gerstenhaber: \emph{On the deformation of rings and 
% algebras II}. Ann. of Math.  84  1966 1--19.

% \bibitem[GER68]{GER3} M. Gerstenhaber: \emph{On the deformation of rings and 
% algebras III}. Ann. of Math. (2)  88  1968 1--34.

% \bibitem[GER74]{GER4} M. Gerstenhaber: \emph{On the deformation of rings and 
% algebras IV}. Ann. of Math. (2)  99  (1974), 257--276.

% \bibitem[GW96]{GER5} M. Gerstenhaber and C. W. Wilkerson: \emph{On the deformation 
% of rings and
% algebras V: Deformation of differential graded algebras}.  Higher homotopy structures in
% topology and mathematical physics (Poughkeepsie, NY, 1996),  89--101, 
% Contemp. Math., 227, Amer. Math. Soc., Providence, RI, 1999.

%\bibitem[HS97]{HiSt} P. J. Hilton and U. A. Stammbach: \emph{A course in homological 
%algebra}. Second edition. Graduate Texts in Mathematics, 4. Springer-Verlag, New York, %1997.

\bibitem{HOC}{G. Hochschild, Cohomology of restricted Lie algebras, \em
Amer. J. Math.,}{\bf 76} (1954), 555--580.

%\bibitem[HS53]{HS} G. Hochschild and J.-P. Serre: \emph{Cohomology of Lie algebras}.
%Ann. of Math. 57 (1953), 591--603.

\bibitem{JAC2}{N. Jacobson, Abstract derivations and Lie algebras, \em
Trans. Amer. Math. Soc.,} {\bf 42} (1937), 206--224.

%\bibitem[JAC80]{JACALG} N. Jacobson: \emph{Basic algebra II}. 
%W. H. Freeman and Co., San Francisco, Calif., 1980.


\bibitem{JAC1}{ N. Jacobson, Classes of restricted Lie algebras of 
characteristic $p$, II, \em Duke Math. Journal,} {\bf 10} (1943), 107-121.  

\bibitem{JACEXE}{ N. Jacobson, \em Exceptional Lie algebras,} 
Lecture Notes in Pure and Applied Mathematics {\bf 1}, Marcel Dekker, New York 1971.


%\bibitem[JAC62]{JACLIE} N. Jacobson: \emph{Lie algebras}. 
%Interscience Tracts in Pure and Applied Mathematics No. 10, New York-London, 1962.

\bibitem{JAN}{ J. C. Jantzen, \em Representations of algebraic groups,} 
Pure and Applied Mathematics, {\bf 131}, Academic Press, Boston, 1987. 


\bibitem{KAC1}{ V. G. Kac, The classification of the simple Lie algebras
over a field with nonzero characteristic (Russian), \em Izv. Akad. Nauk SSSR Ser. Mat.,} 
{\bf 34} (1970), 385--408. English translation: {\em Math. USSR-Izv.,} {\bf 4} (1970), 
391--413.

%\bibitem[KAC71]{KAC2} V. G. Kac: \emph{Global Cartan pseudogroups and simple Lie algebras
%of characteristic $p$} (Russian). Uspekhi Math. Nauk. {\bf 26} (1971). 199--200. 

%\bibitem[KAC74]{KAC3} V. G. Kac: \emph{Description of the filtered Lie algebras with %which graded Lie algebras of Cartan type are associated} (Russian).  
%Izv. Akad. Nauk SSSR Ser. Mat. {\bf 38} (1974), 800--834. 
%English translation: Math. USSR-Izv. {\bf 8} (1974), 801--835.

%\bibitem[KAP71]{KAP} I. Kaplansky: \emph{Lie algebras and locally compact groups}.
%Univ. of Chicago Press, Chicago, 1971.

\bibitem{KS}{ A. I. Kostrikin and I. R. Shafarevich, Cartan's pseudogroups
and the $p$-algebras of Lie (Russian), \em  Dokl. Akad. Nauk SSSR,}{\bf  168}  
(1966), 740--742. English translation: {\em Soviet Math. Dokl.,} {\bf 7} 
(1966), 715--718.

\bibitem{KS2}{ A. I. Kostrikin and I. R. Shafarevich, Graded Lie algebras of
finite characteristic (Russian), \em Izv. Akad. Nauk SSSR Ser. Math.,} {\bf 33} (1969), 
251--322. English translation: {\em Math. USSR-Izv.,} {\bf 3} (1969), 237--304.  

\bibitem{MEL}{ G. M. Melikian, Simple Lie algebras of characteristic $5$
(Russian), \em  Uspekhi Mat. Nauk,}{\bf  35}  (1980), 203--204.

% \bibitem[PS97]{PS1} A. Premet and H. Strade: \emph{Simple Lie algebras of small  
% characteristic. I. Sandwich elements}.  J. Algebra  189  (1997), 419--480.

%\bibitem[PS99]{PS2} A. Premet and H. Strade: \emph{Simple Lie algebras of small %characteristic.
%II. Exceptional roots}.  J. Algebra  216  (1999), 190--301.

%\bibitem[PS01]{PS3} A. Premet and  H. Strade:\emph{Simple Lie algebras of small %characteristic.
%III. The toral rank 2 case}.  J. Algebra  242  (2001), 236--337.

%\bibitem[PS04]{PS4} A. Premet and H. Strade: \emph{Simple Lie algebras of small %characteristic.
% IV. Solvable and classical roots}.  J. Algebra  278  (2004), 766--833.

\bibitem{RUD}{ A. N. Rudakov, Deformations of simple Lie algebras (Russian),
\em Izv. Akad. Nauk SSSR Ser. Mat.,} {\bf 35} (1971), 1113--1119.


%\bibitem[SEL67]{SEL} G. B. Seligman: \emph{Modular Lie algebras}. Ergebnisse der %Mathematik
%und ihrer Grenzgebiete, Band 40. Springer-Verlag, New York, 1967.

\bibitem{SOL1}{ R. Solomon, On finite simple groups and their classification,  
\em Notices Amer. Math. Soc.,}{\bf  42}  (1995), 231--239.

\bibitem{SOL2}{ R. Solomon, A brief history of the classification of the finite 
simple groups, \em  Bull. Amer. Math. Soc.,}{\bf 38}  (2001), 315--352.


%\bibitem[STR89]{STR1} H. Strade: \emph{The Classification of the Simple Modular Lie %Algebras: I.
%Determination of the two-sections}. Ann. of Math. 130 (1989), 643--677.

%\bibitem[STR92]{STR2} H. Strade: \emph{The Classification of the Simple Modular Lie %Algebras:
%II. The Toral Structure}. J. Algebra 151 (1992), 425--475.

%\bibitem[STR91]{STR3} H. Strade: \emph{The Classification of the Simple Modular Lie %Algebras:
%III. Solution of the Classical Case}. Ann. of Math. 133 (1991), 577--604.

%\bibitem[STR93]{STR4} H. Strade: \emph{The Classification of the Simple Modular Lie 
%Algebras: IV. Determining the Associated Graded Algebra}. Ann. of Math. 138 (1993), 
% 1--59.

%\bibitem[STR94]{STR5} H. Strade: \emph{The Classification of the Simple Modular Lie 
% Algebras: V. Algebras with Hamiltonian Two-sections}. Abh. Math. Sem. Univ. Hamburg 64 
%(1994), 167--202.

%\bibitem[STR98]{STR6} H. Strade: \emph{The classification of the simple modular Lie 
% algebras. VI. Solving the final case}.  Trans. Amer. Math. Soc.  350  (1998), 
% 2553--2628.

\bibitem{STR}{ H. Strade, \em Simple Lie algebras over fields of positive 
characteristic I: Structure theory}, De Gruyter Expositions in Mathematics, {\bf 38}, 
Walter de Gruyter, Berlin, 2004.

\bibitem{SW}{ H. Strade and R. L. Wilson, Classification of Simple Lie Algebras over
Algebraically Closed Fields of Prime Characteristic, \em Bull. Amer. Math. Soc.,}
{\bf 24} (1991), 357--362.

\bibitem{VIV1}{ F. Vivani, Infinitesimal deformations of restricted simple 
Lie algebras I, \em J. Algebra} {\bf 320} (2008), 4102--4131.

\bibitem{VIV2}{ F. Viviani, Infinitesimal deformations of restricted simple 
Lie algebras II, \em J. Pure and Applied Algebra} {\bf 213} (2009), 1702-1721.

\bibitem{VIV3}{ F. Viviani, Deformations of the restricted Melikian algebra,
\em Comm. Algebra}  {\bf 37} (2009), 1850--1872. 

%\bibitem{VIV4}{ F. Viviani: Restricted simple Lie algebras and their 
 %infinitesimal deformations, \em Proceedings Conference "From Lie Algebras to Quantum Groups",}
%Centro Internacional de Matematica de Coimbra no. {\bf 28}, 2006 
%(available at arXiv:math.RA/0702755). 

\bibitem{VIV5}{ F. Viviani, Restricted infinitesimal deformations of restricted simple Lie algebras, submitted. }
Preprint available at arXiv:0705.0821.

\bibitem{WAT} { W. C. Watherhouse, \em Introduction to Affine Group Schemes}.
GTM {\bf 66}, Springer-Verlag, 1979. 

\bibitem{WIL}{ R. A. Wilson, \em The finite simple groups}.
GTM {\bf  251}, Springer-Verlag, 2009.

% \bibitem[WEI78]{WEI} B: Ju. Weisfailer: \emph{On the structure of the minimal ideal of 
% some graded Lie algebras of characteristic $p>0$}. J. Algebra {\bf 53} (1978), 
% 344--361. 


%\bibitem[WIL69]{WIL1} R. L. Wilson: \emph{Nonclassical simple Lie algebras}.
%Bull. Amer. Math. Soc. {\bf 75} (1969). 987--991.

% \bibitem[WIL76]{WIL2} R. L. Wilson: \emph{A structural characterization of the simple
% Lie algebras of generalized Cartan type over fields of prime characteristic}.
% J. Algebra {\bf 40} (1976), 418--465.

\end{thebibliography}
\end{document}